\documentclass[11pt, final]{amsart}
\usepackage{amsmath, amsxtra, amssymb, mathdots}
\usepackage{appendix}
\usepackage{verbatim}
\usepackage{multirow}
\usepackage{graphicx}
\usepackage{textcomp}

\usepackage[svgnames]{xcolor}
\usepackage{tikz}
\usepgflibrary{decorations.text}
\usetikzlibrary{decorations.text}
\pgfdeclarelayer{edgelayer}
\pgfdeclarelayer{nodelayer}
\pgfsetlayers{edgelayer,nodelayer,main}

\tikzstyle{none}=[inner sep=0pt]
\definecolor{hexcolor0xfefdfd}{rgb}{0.996,0.992,0.992}

\usepackage[margin=1.3in]{geometry}
\usepackage{graphicx}
\usepackage[cmtip,all]{xy}
\usepackage{xy}

\usepackage[notcite, notref]{showkeys}
\usepackage[colorlinks, linkcolor=blue, citecolor=red, urlcolor=black]{hyperref}

\theoremstyle{plain}
\newtheorem{main-theorem}{Main Theorem}
\newtheorem{theorem}[equation]{Theorem}

\newtheorem*{claim*}{Claim}

\theoremstyle{definition}

\newtheorem{remark}[equation]{Remark}

\numberwithin{equation}{section}

\title{Nontransitive Identities: A New Method for Creating Nontransitive Dice}
\author{Logan Danielson}

\begin{document}
\maketitle

\begin{abstract}
In this article, a new method for characterizing nontransitive dice is described. This new method is then used to describe the \textquotesingle \textquotesingle Nontransitive Identities\textquotesingle \textquotesingle \ (NI) that are possible for 3 dice with 3, 4 and 5 sides each as well as for 5 dice with 3 sides each. Next, we will discuss how these NI can be used to create NI that involve more dice and/or die sides. From there, the 3 dice Nontransitive Identities (NI) will be used to produce 5 dice NI. We will begin this by describing more nomenclature that becomes necessary for characterizing sets of NI with multiple win chains. Then, we will describe a method for creating 5 dice NI and 7 dice NI from our 3 dice NI. We will conclude with a discussion on the patterns that exist in the frequency of nontransitive identities and how these patterns correlate to NI morphisms. Our ancillary files will include our python scripts as well as all of the sets of Nontransitive identities we have discussed.
\end{abstract}

\section{Introduction}

Nontransitive dice have been a subject of fascination for many mathematicians and lay-people, like myself, ever since the Efron Dice were mentioned in the Scientific American magazine. [1] My first introduction to them was a Numberphile video titled \textquotesingle \textquotesingle The Most Powerful Dice.\textquotesingle \textquotesingle \  Nontransitive dice are dice with non-standard face values. The dice are used to play a game where the winner is person whose die rolls the highest face value. The dice are said to be \textquotesingle \textquotesingle Nontransitive\textquotesingle \textquotesingle \  because for any die in the set, there is always another die in the set that will more frequently roll a higher value. In this way, the nontransitive dice captivate us by violating our seemingly innate assumption of \textquotesingle \textquotesingle this is better than that\textquotesingle \textquotesingle \ relationships being transitive. Rock-paper-scissors is another example of a game that is nontransitive because there is no choice that is the best choice, if don\textquotesingle t know the opponent\textquotesingle s choice— or, stated differently: the first choice is the worst choice. 

In this paper, we will only concern ourselves with \textquotesingle \textquotesingle perfectly\textquotesingle \textquotesingle \ nontransitive dice. That is to say, every die in a nontransitive set of dice will win and lose to the same number of dice. Pictorially, our nontransitive dice will be represented complete graphs, whose vertices are labelled with letters that represent each die involved in the nontransitive set of dice. These vertices are connected by an arrow that designates which die \textquotesingle \textquotesingle beats \textquotesingle \textquotesingle \ the other by more frequently rolling a higher value. This arrow is directed towards the die that loses. For three dice, we have: 

\begin{figure}[h!]
\includegraphics[width=\textwidth]{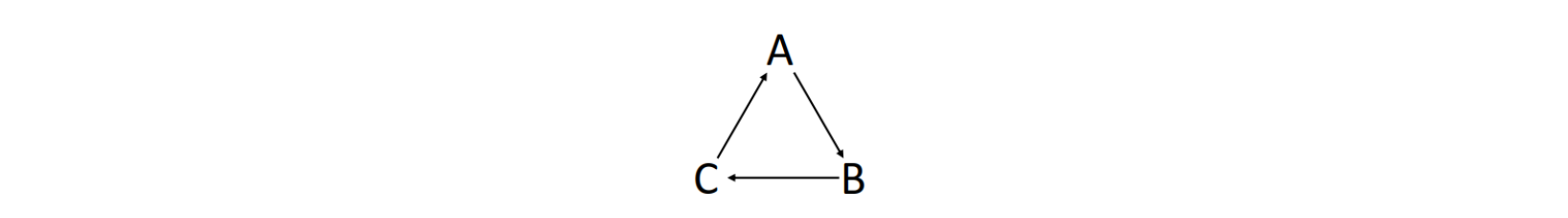}
\end{figure}

In this scenario, Dice A \textquotesingle \textquotesingle beats\textquotesingle \textquotesingle , i.e. more frequently rolls a higher value than, Dice B. Dice B beats Dice C and Dice C beats Dice A. If the dice were \textquotesingle \textquotesingle transitive \textquotesingle \textquotesingle \ dice, then this chain would be broken and Dice A would beat Dice C. The fact that there is no ultimate winner to the chain is what makes the set of dice \textquotesingle \textquotesingle nontransitive.\textquotesingle \textquotesingle  

Part of the beauty of the math involved with nontransitive dice is that understanding it only truly requires an understanding of what value is greater and the ability to tally wins. For this reason, I will endeavor to refrain from using complicated notation in this article and I intend to restrict my use of the word \textquotesingle \textquotesingle probability \textquotesingle \textquotesingle \ to just this sentence. 

Previously, many sets of nontransitive dice have been created. Many of these sets were created via trial and error and others were generated via a formula to produce nontransitive dice of a particular pattern. [2] This paper hopes to elucidate the patterns that underlie nontransitive sets of dice and to define these patterns in that same basic sense of what is \textquotesingle \textquotesingle greater than\textquotesingle \textquotesingle \ something else. 

After describing these patterns, that I\textquotesingle ll refer to as \textquotesingle \textquotesingle Nontransitive Identities\textquotesingle \textquotesingle ,  I\textquotesingle ll describe several operations that can be performed with nontransitive identities to generate sets of perfectly nontransitive dice with more sides as well as sets of dice involving more dice and more sides. 

Later in this article, we will use a pictorial representation of the non-transitive dice to explain how we can create 5 dice NI from 3 dice NI. As the number of dice in a NI increases, the need for more nomenclature regarding NI becomes apparent and for that reason, we will define more nomenclature for describing NIs. This nomenclature will be used to help characterize the 5 dice nontransitive identities generated by overlapping 3 dice NI. Additionally, a proof will be shown for why the change in the secondary win chain for 5 dice identities does not result in a unique nontransitive identity.

We\textquotesingle ll conclude with a discussion on the patterns that exist in the frequency of nontransitive identities. To imply that the nontransitive identities have a frequency is to say that the NI can be mapped to some sort of one dimensional spectrum. In effect, we forced the NI along a spectrum and, in this article, we\textquotesingle ll discuss how this was done and the nature of patterns that were found. From the outset, we acknowledge that this method is biased to detecting patterns between NI that exist along the right end of an NI and as such, this method is definitely not perfect. It is hoped that the perception of patterns that can be seen by this method will inspire others to characterize the higher-dimensional categories of nontransitive identities.

\section{Methodology}

The math in this paper was performed with a python script. The first iteration of the script tried every possible combination of numbers on 3 dice with 3 sides, played them against each other and returned the combinations that made the dice nontransitive. These sets of dice were then used to characterize the list of identities. The number of combinations of sets of 3 dice with 4 sides was too large to process on my personal computer in a reasonable timeframe. To address this problem, I rewrote the script to iterate and test \textquotesingle \textquotesingle viable\textquotesingle \textquotesingle \ identities instead of directly iterating the numbers on the dice faces themselves. An identity was defined as \textquotesingle \textquotesingle viable\textquotesingle \textquotesingle \  if it coded for the appropriate number of die faces for each die. The 3rd version of my script used a significantly faster method to idenitify the viable identities to be tested.

\section{The Nontransitive identities}		

These identities, when \textquotesingle \textquotesingle solved\textquotesingle \textquotesingle \ produce a set of nontransitive dice. The process of \textquotesingle \textquotesingle solving\textquotesingle \textquotesingle \ the identity involves reading the identity from left to right and assigning values to a face of the die indicated by the letter, in accordance with the pattern established by the nontransitive identity. For example, this is a nontransitive identity for a set of 3 dice with 3 sides each:
\begin{figure}[h!]
\includegraphics[width=\textwidth]{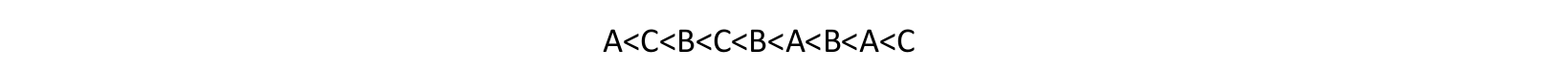}
\end{figure}

\noindent It can be solved by assigning the values to the dice faces as follows:

Dice A\textquotesingle s sides: 1, 6, 8

Dice B\textquotesingle s sides: 3, 5, 7

Dice C\textquotesingle s sides: 2, 4, 9

\

\noindent Alternatively, it could be solved by assigning the values of:

Dice A\textquotesingle s sides: 1, 10, 14

Dice B\textquotesingle s sides: 7, 9, 13

Dice C\textquotesingle s sides: 4, 8, 22

\

\noindent Indeed, there are an infinite number of ways that the identity could be solved, and yet, with the nomenclature of the identity, we can characterize that infinite number of sets of nontransitive dice with just one nontransitive identity.

These identities are written such that the sequence of who beats who, hereby referred to as the \textquotesingle \textquotesingle win chain\textquotesingle \textquotesingle,  proceeds alphabetically and the last letter loops back to beat the first. As the win chain has no true start and no true end, it was arbitrarily decided for Die A to be the die with the lowest value die face because, for any set of nontransitive dice, one of the dice would have to have the face with the lowest value. In the event of multiple dice sharing the lowest value, then the subset of dice that shared the lowest value is checked for the next lowest value and Die A is defined as being that die. If this next lowest value is also shared, then this process of identifying a subset of dice and checking for their next lowest value reiterates until next lowest value is held by the die that is die A alone. 

It also must be noted that these identities were written to minimize the value of the highest value die face. This is important because any portion within the identity of the form \textquotesingle \textquotesingle \textless N=N\textquotesingle \textquotesingle ,  where \textquotesingle \textquotesingle N=N\textquotesingle \textquotesingle \  refers to two faces of die N, could be replaced with \textquotesingle \textquotesingle \textless N\textless N\textquotesingle \textquotesingle \ without affecting the nontransitive characteristic of the identity. While the change does not affect the nontransitive characteristic, it does enable some nontransitive identities to be expanded in multiple ways with some of the expansion techniques that will be discussed later. The identities are also written such that dice faces of equivalent value are listed in alphabetical order. That is to say, if a portion of an identity read \textquotesingle \textquotesingle \textless A=B=C\textquotesingle \textquotesingle \ and it was replaced with \textquotesingle \textquotesingle \textless B=A=C\textquotesingle \textquotesingle ,  then the modified identity would produce the same infinite combination of sets of dice as the original identity. To prevent this duplicity, nontransitive identities with dice faces of equivalent value are only listed with these equivalent faces in alphabetical order.

The nontransitive identities (NI) could potentially be more specifically referred to as \textquotesingle \textquotesingle nontransitive identity inequalities\textquotesingle \textquotesingle \ but the NI are not always completely \textquotesingle \textquotesingle inequal\textquotesingle \textquotesingle \ and the addition of the word \textquotesingle \textquotesingle inequality\textquotesingle \textquotesingle \ doesn\textquotesingle t add meaning to the concept conveyed by the shorter name, NI. For this reason, I\textquotesingle ll use the shorter name, NI. Also, \textquotesingle \textquotesingle nontransitive identities\textquotesingle \textquotesingle \  sounds like something a knight would say.

I\textquotesingle d also imagine that some people would prefer to see these identities written with numbered subscripts. This makes a certain amount of sense because algebra tells us that, when we \textquotesingle \textquotesingle solve\textquotesingle \textquotesingle \ an equation, that one variable will receive one value. As such when we solve a nontransitive identity and in the process assign one letter multiple values, it just feels wrong for a bit and we could escape that feeling by adding subscripts to the letters. My reason for not doing this is two-fold. First, this is not elementary algebra. Second, numbers are ordinal and if you subscript with numbers, then you run the risk creating a subtle implied distinction between the die faces, a distinction based on our understanding of the ordinal nature of numbers. Given that no ordinal distinction between the die faces exists in the context of nontransitive identities and nontransitive dice, our nomenclature should support that lack of a distinction. That said, you will see the NI with numbered subscripts twice in this paper. This was done to more clearly explain how the NI was expanded in these two cases.

\section{Nontransitive identities for 3 dice with 3 sides}		

\begin{figure}[h!]
\includegraphics[width=\textwidth]{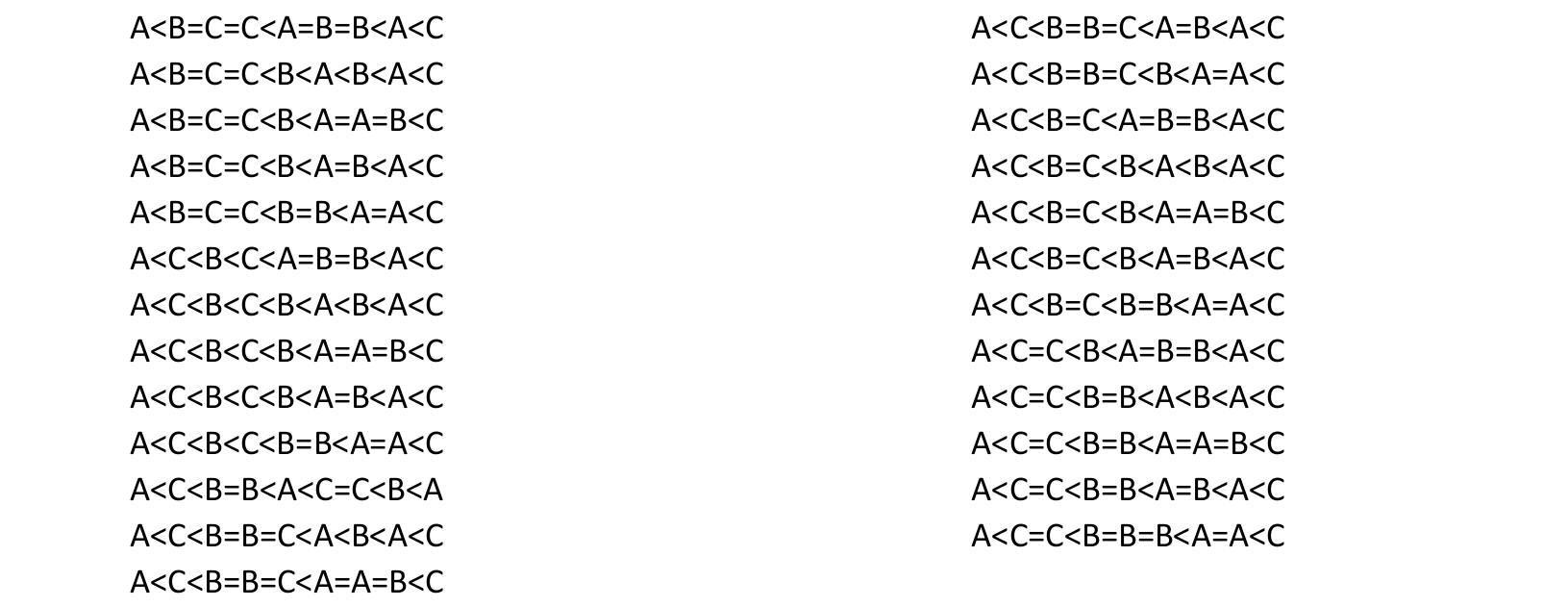}
\end{figure}

\begin{remark} 
In the ancillary files for this article, the complete sets of Nontransitive Identities have been identified for 3 dice with 3, 4 and 5 sides each as well as for 5 dice with 3 sides each. The \textquotesingle \textquotesingle Alphabetical Dupe\textquotesingle \textquotesingle \ versions contain the duplicitous "\textless N\textless N" NI but do not allow the equivalent dice faces to be out of alphabetical order. The \textquotesingle \textquotesingle Irreducible \textquotesingle \textquotesingle \ versions do not contain any duplicitous NI.  
\end{remark}

\pagebreak

\section{Eponymous Nontransitive Identities}		

\begin{figure}[h!]
\includegraphics[width=\textwidth]{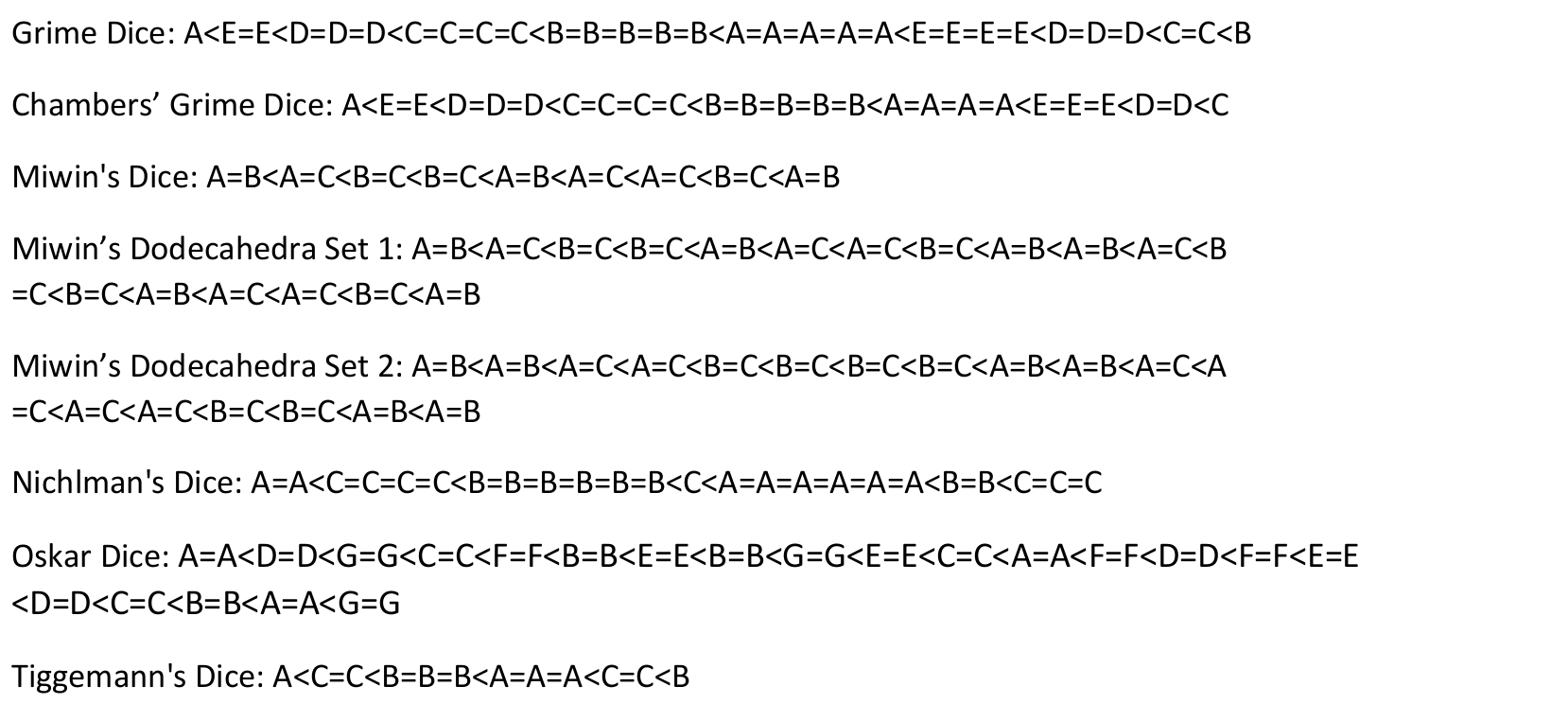}
\end{figure}

\section{Methods for Expansion of Nontransitive Identities}

For this section, I\textquotesingle ll state the method as a theorem and then either provide a logical proof of the method or an example of the method\textquotesingle s use. In the later case, I will have tested the method thoroughly via python scripts and the details of what is meant by \textquotesingle \textquotesingle thorough\textquotesingle \textquotesingle \ will be provided in the text.  

\begin{theorem}
\textquotesingle \textquotesingle Add Zero\textquotesingle \textquotesingle :  Nontransitive identities will remain nontransitive when the end of the identity is appended with \textquotesingle \textquotesingle \textless A=B … =N\textquotesingle \textquotesingle ,  where N is the last die of the nontransitive identity. 
\end{theorem}

This theorem implies that \textquotesingle \textquotesingle A\textless C\textless B\textless C\textless B\textless A\textless B\textless A\textless C\textquotesingle \textquotesingle \ + \textquotesingle \textquotesingle \textless A=B=C\textquotesingle \textquotesingle \ produces nontransitive identity, and it does given that \textquotesingle \textquotesingle A\textless C\textless B\textless C\textless B\textless A\textless B\textless A\textless C\textless A=B=C\textquotesingle \textquotesingle \  is in our set of nontransitive identities for 3 dice with 4 sides. Logically, this expansion makes sense because it results in the number of wins for each die increasing by the same amount.

\begin{theorem}
\textquotesingle \textquotesingle Multiply by One\textquotesingle \textquotesingle :  Nontransitive identities will remain nontransitive when each die face of the identity, designated by the letter \textquotesingle \textquotesingle N\textquotesingle \textquotesingle ,  is replaced by \textquotesingle \textquotesingle N=N\textquotesingle \textquotesingle \ or \textquotesingle \textquotesingle N\textless N \textquotesingle \textquotesingle.
\end{theorem}

This theorem implies that \textquotesingle \textquotesingle A=A\textless C=C\textless B=B\textless C=C\textless B=B\textless A=A\textless B=B\textless A=A\textless C=C\textquotesingle \textquotesingle \ is a nontransitive identity, given that \textquotesingle \textquotesingle A\textless C\textless B\textless C\textless B\textless A\textless B\textless A\textless C\textquotesingle \textquotesingle \ is a nontransitive identity. This expansion makes sense because it results in the number of wins for each die multiplying by the same amount.

\pagebreak

\begin{theorem}
\textquotesingle \textquotesingle Identity Addition\textquotesingle \textquotesingle :  Nontransitive identities will remain nontransitive when the end of the identity is appended with another Nontransitive identity of the same number of dice and die sides. Multiple nontransitive identities can be appended at once and between each can exist a \textquotesingle \textquotesingle =\textquotesingle \textquotesingle \ or a \textquotesingle \textquotesingle \textless \textquotesingle \textquotesingle . 
\end{theorem}

This theorem was checked for every combination of nontransitive identities (NI) that could be created from the 3 dice and 3 side NI list, using either two or three nontransitive identities and either a \textquotesingle \textquotesingle =\textquotesingle \textquotesingle \ or a \textquotesingle \textquotesingle \textless \textquotesingle \textquotesingle \  between each identity. All combinations produced new nontransitive identities. This theorem was also checked for 28 million NI created using all combinations of \textquotesingle \textquotesingle \textless \textquotesingle \textquotesingle \ and  \textquotesingle \textquotesingle = \textquotesingle \textquotesingle \ between 9 NI.  Again, all combinations produced new nontransitive identities.

\subsection{Identity Addition example:}
\	
\begin{figure}[h!]
\includegraphics[width=\textwidth]{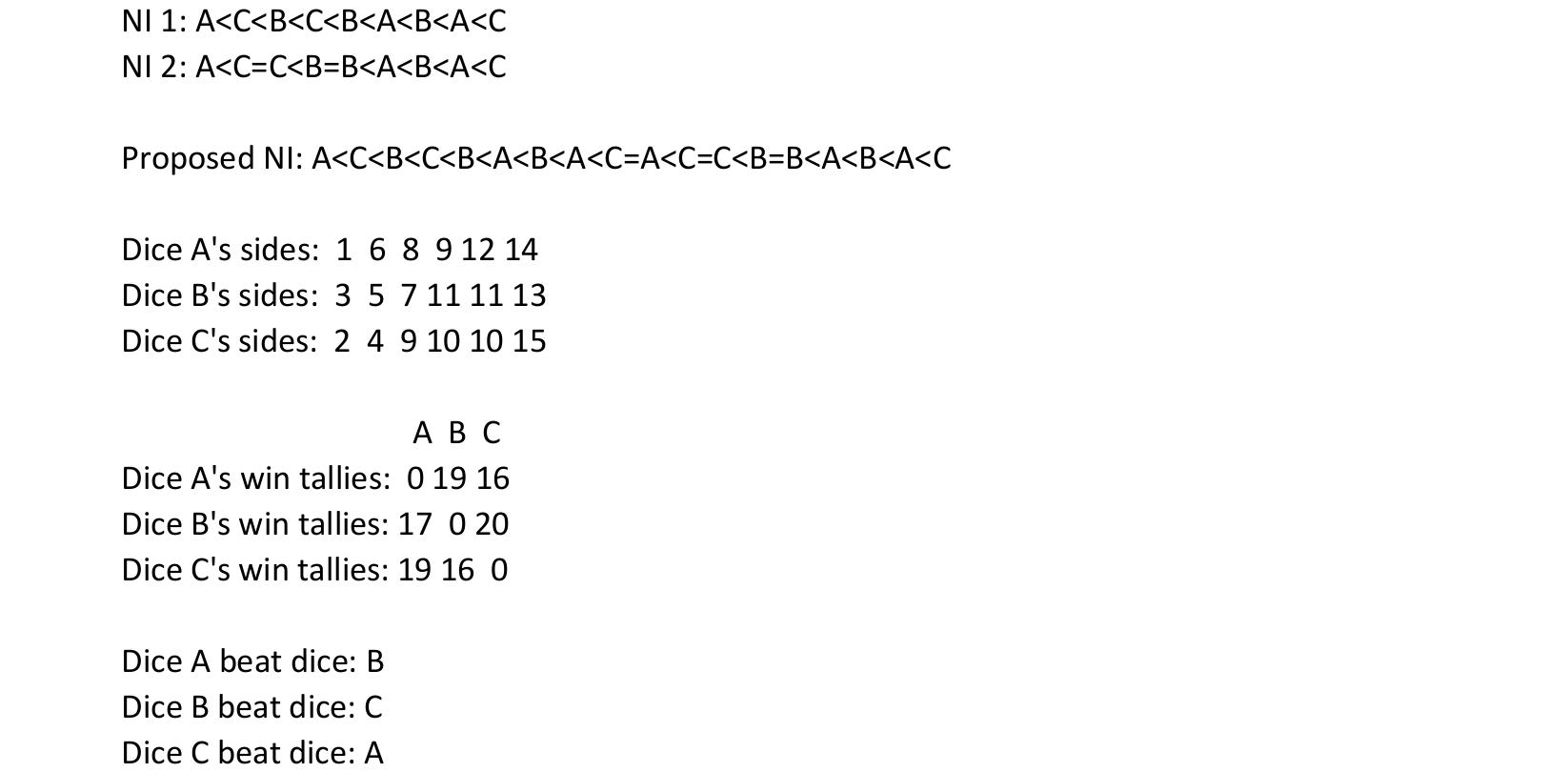}
\end{figure}

\newpage

\begin{theorem}
\textquotesingle \textquotesingle Identity Nesting\textquotesingle \textquotesingle :  Nontransitive identities will remain nontransitive when each die face of the \textquotesingle \textquotesingle base\textquotesingle \textquotesingle \  nontransitive identity, designated by the letter \textquotesingle \textquotesingle N\textquotesingle \textquotesingle ,  is replaced by a \textquotesingle \textquotesingle substituting\textquotesingle \textquotesingle \ nontransitive identity, where each of the substituting nontransitive identities have the same number of dice and dice faces.
\end{theorem}

This method provides us with two manners by which we can expand a given nontransitive identity. The manners differ in whether they prioritize increasing the number of dice faces or if they increase the number of both the dice and dice faces. The identity nesting methods also features certain special cases that will impose some rules on how a given base identity can be expanded. These topics will be all be briefly introduced and then discussed in more detail sequentially.

Theorem 6.4a \textquotesingle \textquotesingle Identity Face Exponentiation\textquotesingle \textquotesingle :  Nontransitive identities will remain nontransitive when each die face of the \textquotesingle \textquotesingle base\textquotesingle \textquotesingle \ nontransitive identity, designated by the letter \textquotesingle \textquotesingle N\textquotesingle \textquotesingle ,  is replaced by a \textquotesingle \textquotesingle substituting\textquotesingle \textquotesingle \ nontransitive identity, where each of the substituting nontransitive identities have the same number of dice and dice faces.

This produces a nontransitive identity with the number of dice equal to the involved in the substituting nontransitive identities and with a number of die sides equal to the number of die sides of the base nontransitive identity raised to the power of the number of die sides of the substituting nontransitive identities.

Theorem 6.4b \textquotesingle \textquotesingle Identity Dice Multiplication\textquotesingle \textquotesingle :  Nontransitive identities will remain nontransitive when each die face of the identity, designated by the letter \textquotesingle \textquotesingle N\textquotesingle \textquotesingle ,  is replaced by a nontransitive identity, where each of the substituting nontransitive identities have the same number of dice and dice faces. And, where the substituting nontransitive identities for all dice of the base nontransitive identity, except that of Die A, are reassigned previously unused letters.

This produces a nontransitive identity with a number of dice and die sides equal to the product of the number involved in the base nontransitive identity and substituting nontransitive identities, respectively.

Theorem 6.4\textquotesingle s special cases: When a portion of the base nontransitive identity is of the form of \textquotesingle \textquotesingle \textless N=N\textquotesingle \textquotesingle \ or \textquotesingle \textquotesingle \textless M=N\textquotesingle \textquotesingle \ where N and M refer to dice faces of different dies, then the expansion of this base nontransitive identity must maintain this special relationship. This is achieved by imposing that both of these die faces be expanded with the same substituting nontransitive identity and in a manner similar to the Theorem 2\textquotesingle s \textquotesingle \textquotesingle Multiply by One\textquotesingle \textquotesingle \ operation. 

\newpage

\subsection{Identity Face Exponentiation example:}
\		
\begin{figure}[h!]
\includegraphics[width=\textwidth]{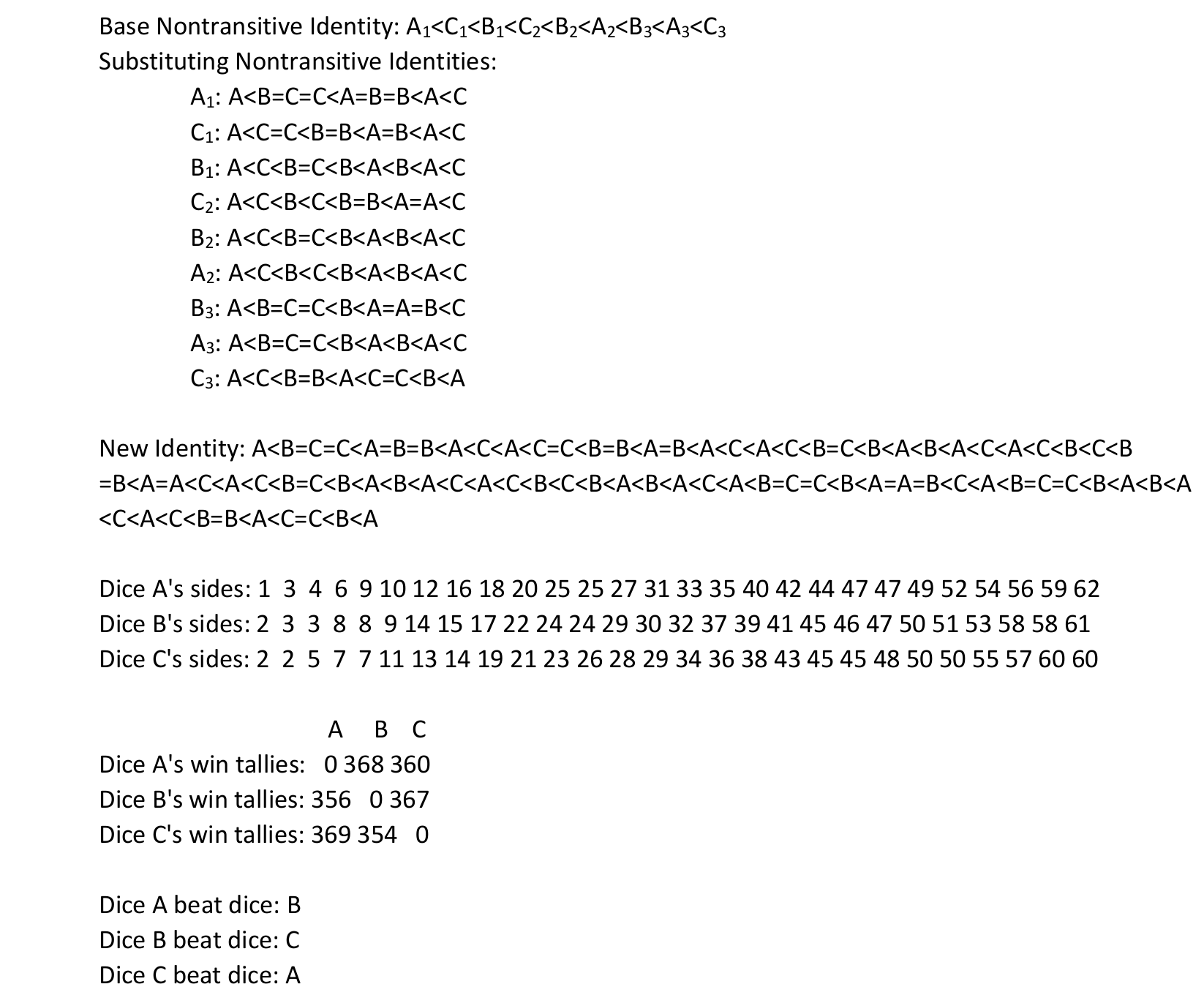}
\end{figure}

Theorem 6.4a\textquotesingle s Identity Face Exponentiation operation was checked for every combination of nontransitive identities (NI) that could be created from the 3 dice and 3 side NI list, except for combinations whose base nontransitive identities fit the \textquotesingle \textquotesingle \textless M=N\textquotesingle \textquotesingle \  special case. To ease the computation load, each die face of a particular die was expanded with the same NI for all faces of that die. All combinations produced new nontransitive identities.

\newpage

\subsection{Identity Dice Multiplication example:} 
\
\begin{figure}[h!]
\includegraphics[width=\textwidth]{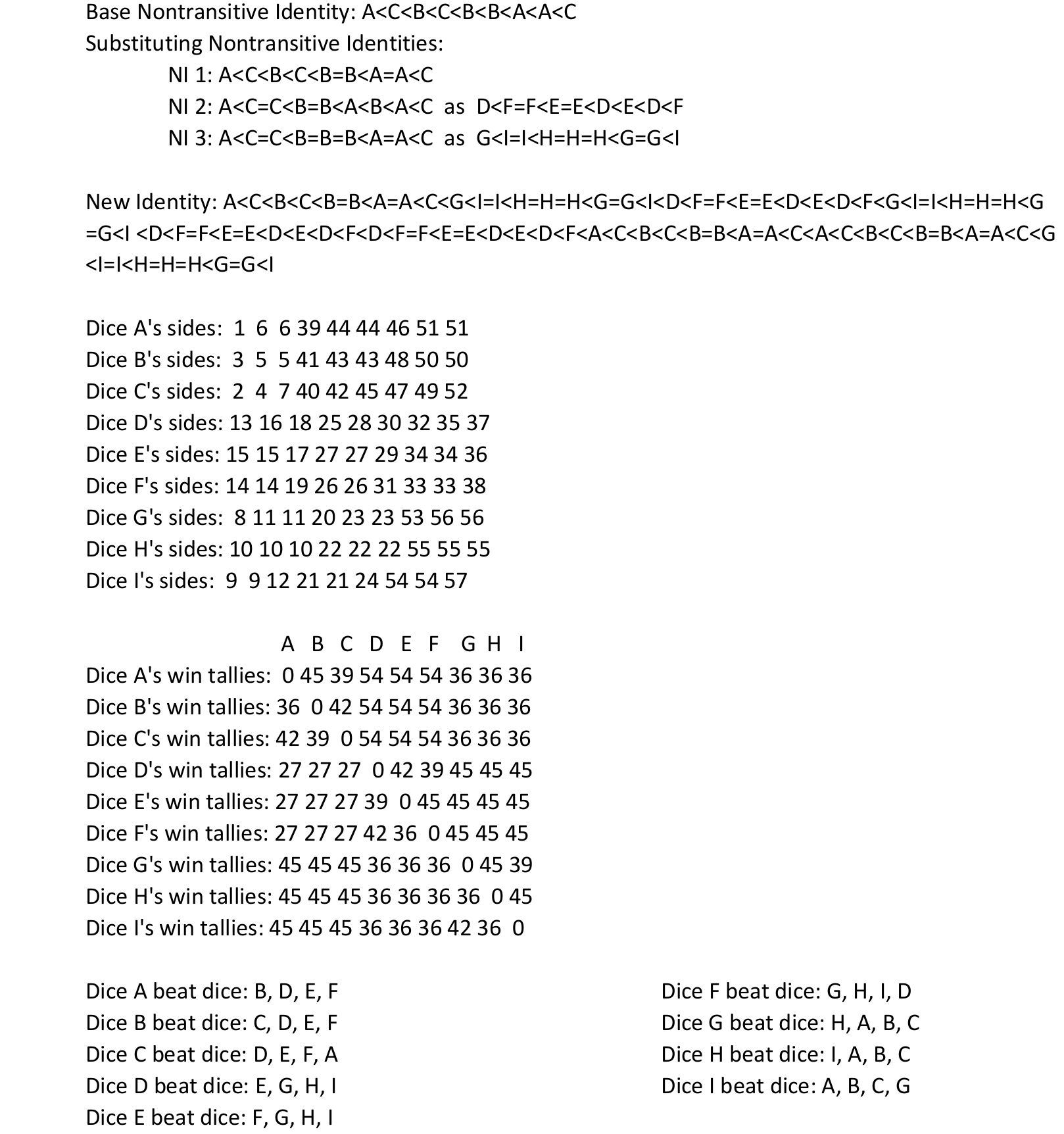}
\end{figure}

Theorem 6.4b\textquotesingle s Identity Dice Multiplication operation was checked for every combination of nontransitive identities (NI) that could be created from the 3 dice and 3 side NI list, except for combinations whose base nontransitive identities fit the \textquotesingle \textquotesingle \textless M=N\textquotesingle \textquotesingle \ special case. To ease the computation load, each die face of a particular die was expanded with the same NI. All combinations produced new nontransitive identities.

You may have noticed that the base NI used in the previous example was not included in the list of NIs for 3 dice and 3 sides. This is because it was not listed. That base NI was a duplicative form of the one used in the following example. Please note, the previous example was one method of how to expand a base nontransitive identity that has a portion of the form \textquotesingle \textquotesingle \textless N=N\textquotesingle \textquotesingle ,  that is to say, the \textquotesingle \textquotesingle \textless N=N\textquotesingle \textquotesingle \  is replaced with \textquotesingle \textquotesingle \textless N\textless N\textquotesingle \textquotesingle \  and expanded as any other die face would be. The following example shows the other method for expanding this type of special case.

\subsection{Theorem 6.4 special case \textquotesingle \textquotesingle \textless N=N\textquotesingle \textquotesingle \ example:}		
\
\begin{figure}[h!]
\includegraphics[width=\textwidth]{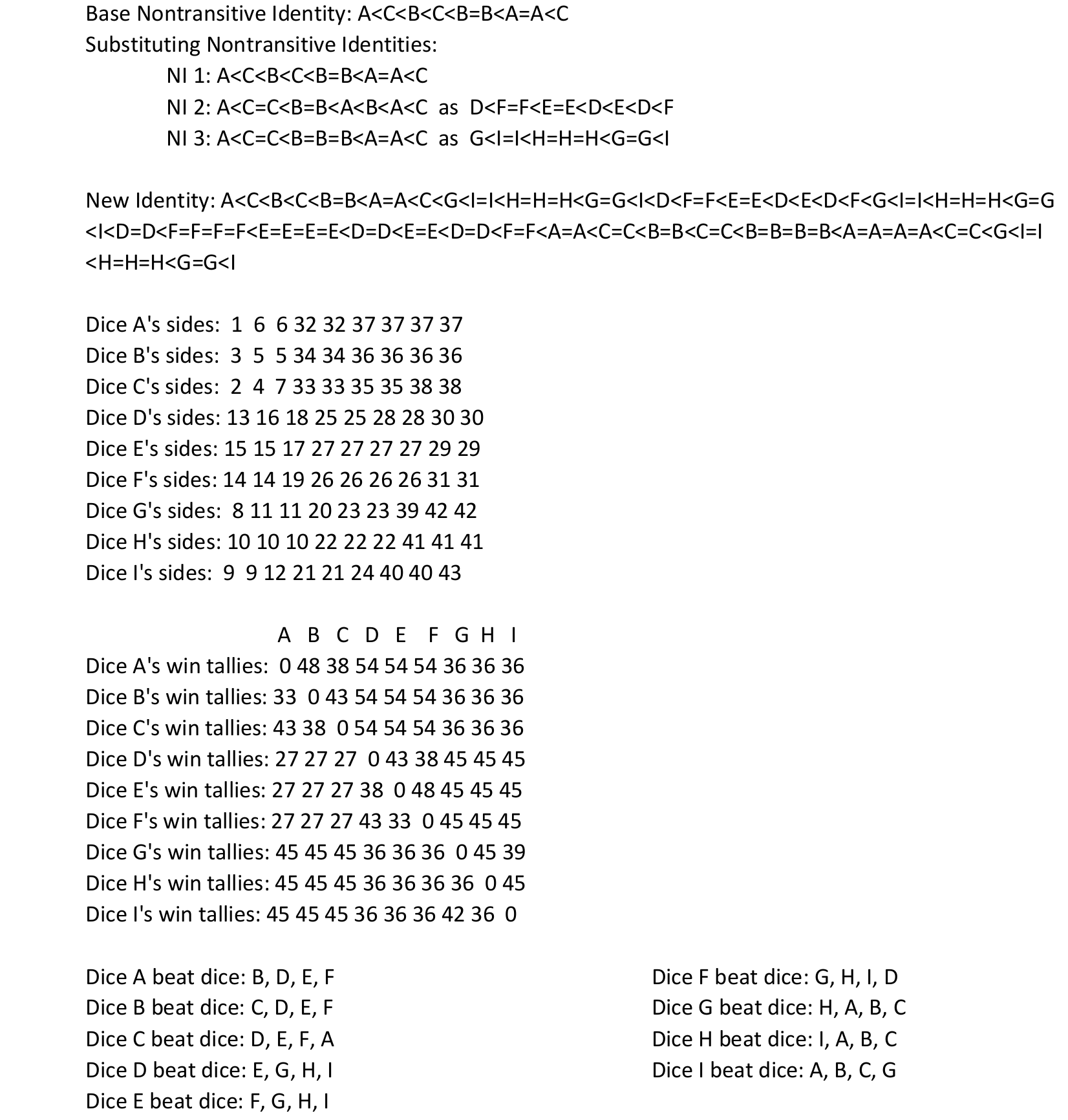}
\end{figure}

\newpage

\subsection{Theorem 6.4 special case \textquotesingle \textquotesingle \textless M=N\textquotesingle \textquotesingle \ example:}		
\
\begin{figure}[h!]
\includegraphics[width=\textwidth]{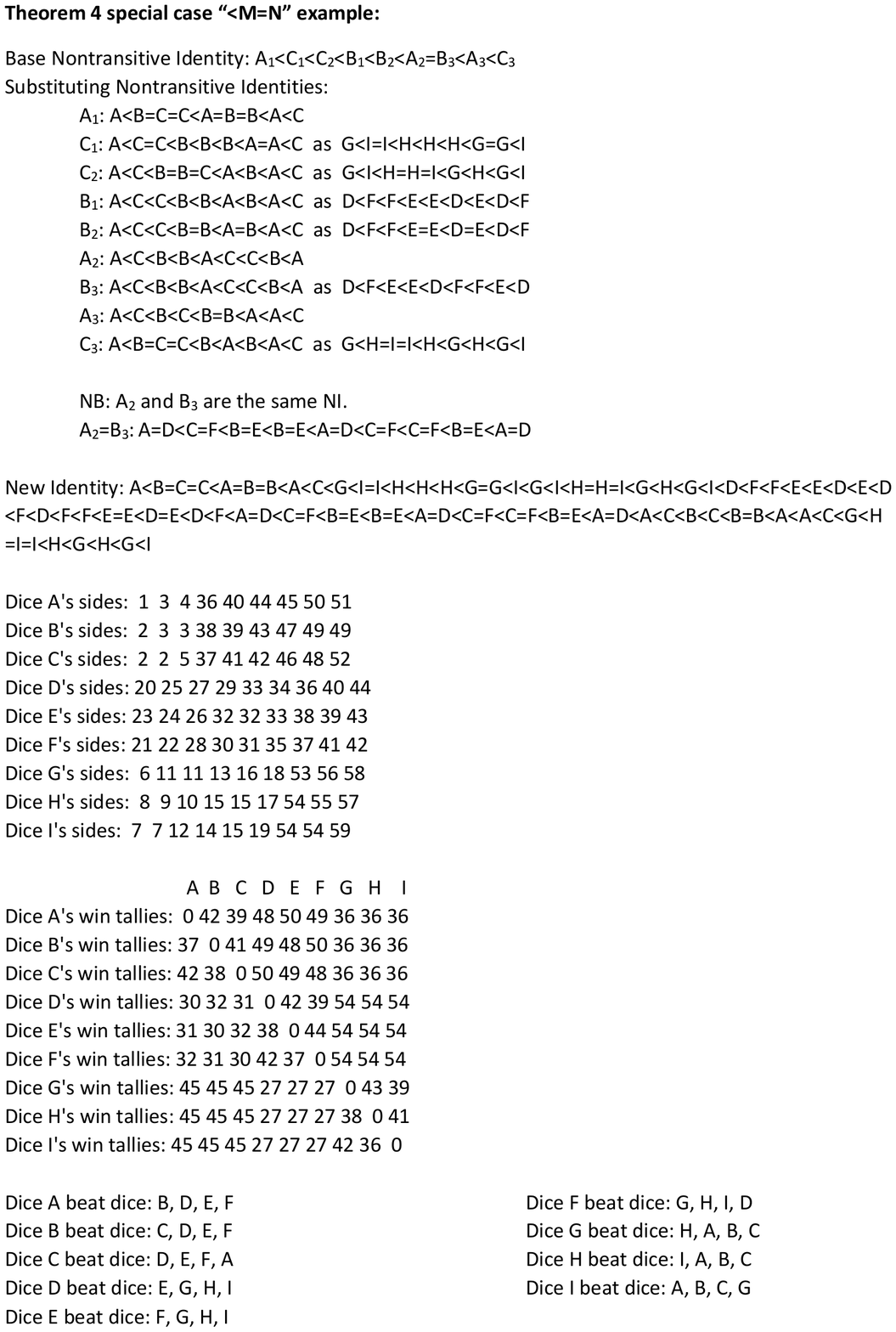}
\end{figure}

\subsection{Pictorial Representations of Expansions:} 		
The nontransitive identities produced by the Identity Dice Multiplication method do follow what we would predict from their pictorial representations. 
\
\begin{figure}[h!]
\includegraphics[width=\textwidth]{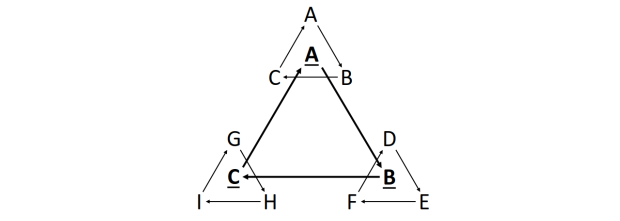}
\end{figure}

In the above image, the underlined letters represent the dice of the base NI and the other letters represent the dice that are present in the NI that is produced by the substituting NI being nested in the base NI. It follows from this picture that all dice created by nesting into the Base NI\textquotesingle s die A (i.e. Dice A, B and C) should beat that of dice of B (i.e. Dice D, E and F), and that all dice of B should beat the dice of C (i.e. Dice G, H and I), and that all the dice of C should beat the dice of A. This is exactly what happened in our examples and it can be represented as an enneagram as follows.

\begin{figure}[h!]
\includegraphics[width=\textwidth]{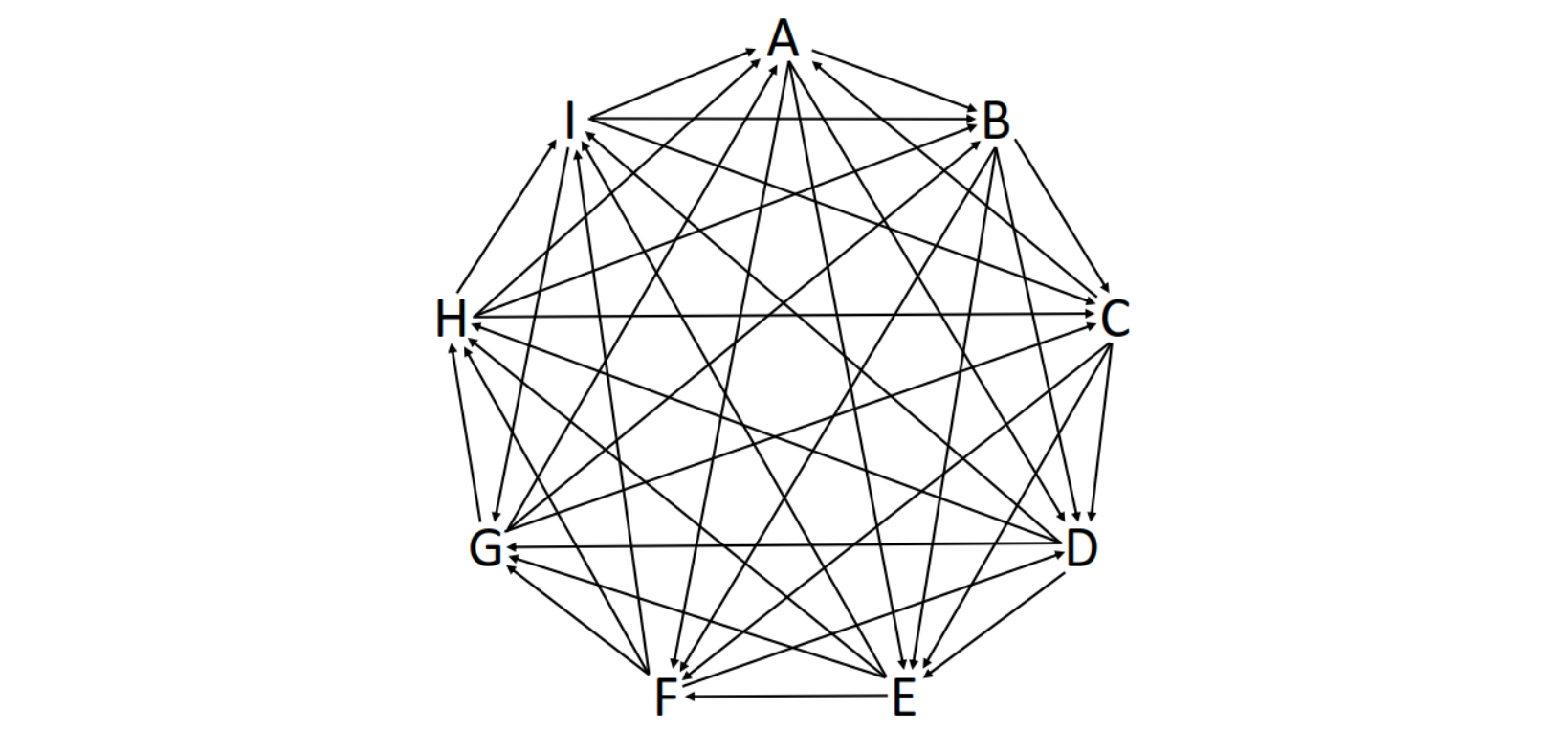}
\end{figure}

\pagebreak

\section{Additional Nontransitive Identity nomenclature}
First, please consider the following pictorial representations of 5 dice identity win chains.
\begin{figure}[h!]
\includegraphics[width=\textwidth]{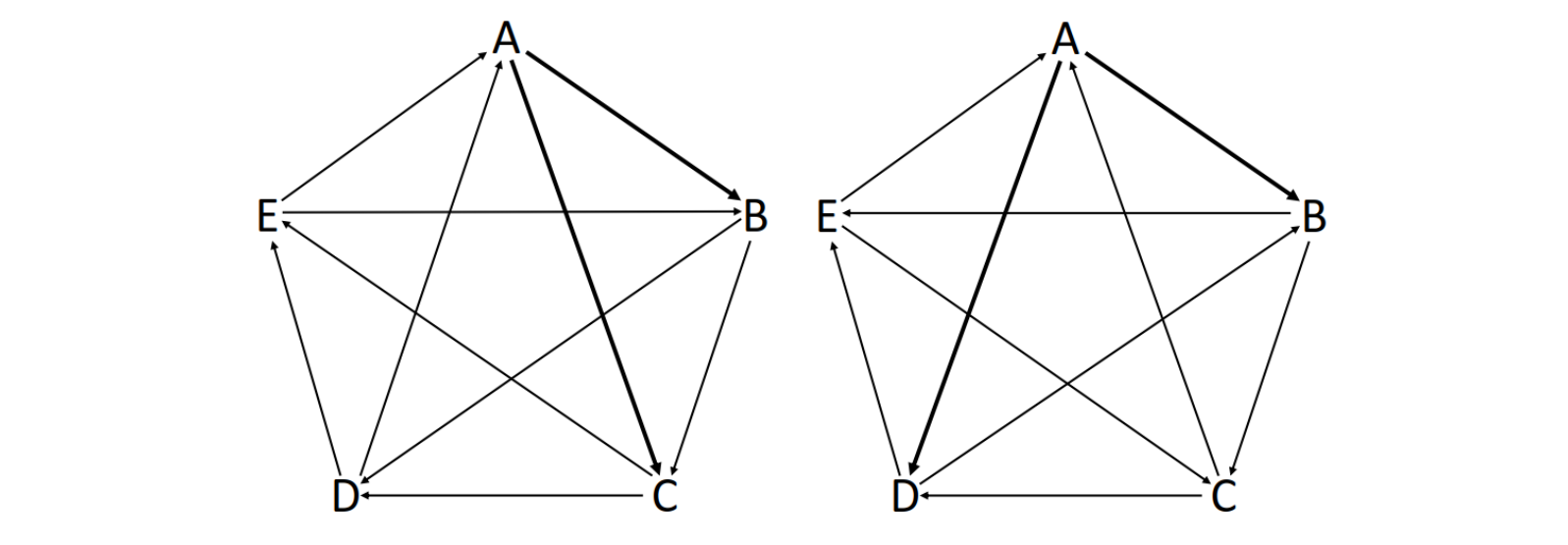}
\end{figure}

In these pictures, we see that the primary win chain proceeds alphabetically, but the secondary win chain either directly follows the primary win chain or skips a letter prior to starting the next win chain. For 5 dice NI, this difference does not matter and this will be proven later in this article. For now, let\textquotesingle s assume that the two complete graphs represent different NI so that we can start to delineate our new nomenclature with easier examples.

\

The Grime dice follows the win chain pattern of the pentagram on the left. The pattern of the win chain for this identity could be described as follows:

\begin{figure}[h!]
\includegraphics[width=\textwidth]{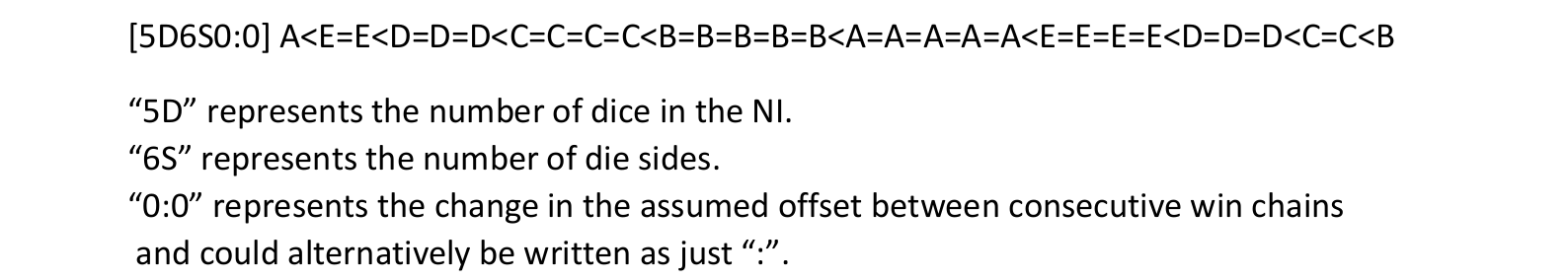}
\end{figure}

\noindent By \textquotesingle \textquotesingle assumed offset\textquotesingle \textquotesingle ,  I mean to say that our NI nomenclature assumes a particular pattern between successive win chains, i.e. that the primary win chain proceeds alphabetically, in a clockwise manner, and that subsequent win chains result in a die beating the following die along the previous win chain. 

\

The win chain pattern described by the pentagram on the right could be written as:

\begin{figure}[h!]
\includegraphics[width=\textwidth]{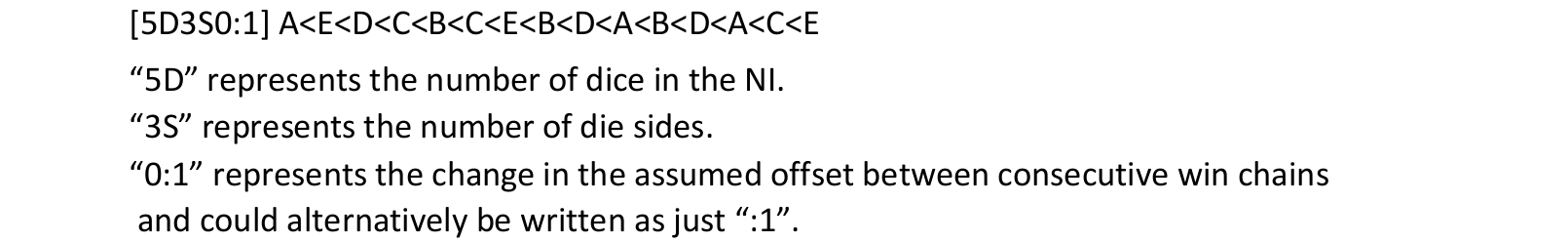}
\end{figure}

\noindent The first zero in \textquotesingle \textquotesingle 0:0\textquotesingle \textquotesingle \ implies that the first win chain proceeds alphabetically. For characterized NI, this should always be the case. But, for theoretical NI sets produced by trying to add or remove dice from other sets of NI, the letters of the primary win chain can be placed after the \textquotesingle \textquotesingle ]\textquotesingle \textquotesingle \  to define the order of the primary win chain, as if they were in alphabetical order, with the die of the lowest value face listed first. If this nomenclature was applied to our assumed primary win chain, then our \textquotesingle \textquotesingle ]\textquotesingle \textquotesingle \ would follow \textquotesingle \textquotesingle ] ABCDE\textquotesingle \textquotesingle . The second zero in \textquotesingle \textquotesingle 0:0\textquotesingle \textquotesingle \  implies that the second win chain follows immediately after the first. The one in \textquotesingle \textquotesingle 0:1\textquotesingle \textquotesingle \ implies that the second win chain skips a die, or \textquotesingle \textquotesingle steps\textquotesingle \textquotesingle \  over a die, i.e. there\textquotesingle s a \textquotesingle \textquotesingle step\textquotesingle \textquotesingle \ in the win chain, prior to starting the next win chain. 

\

The Oskar dice produce an NI that can be more fully described as:

\begin{figure}[h!]
\includegraphics[width=\textwidth]{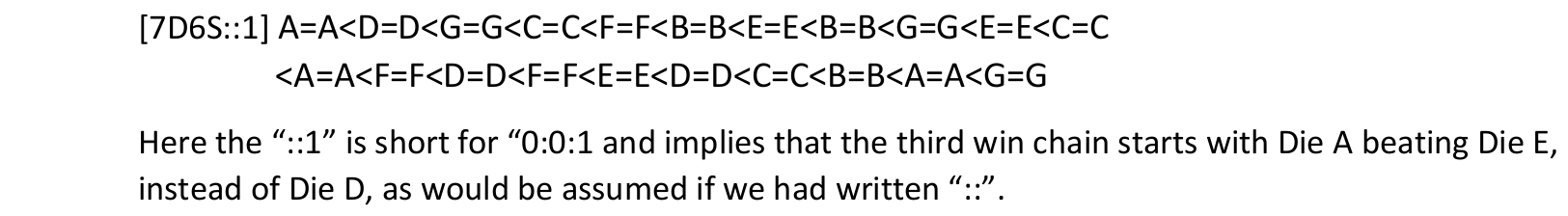}
\end{figure}
	
The 9 Dice NI created by the Identity Dice Multiplication method described in my previous article [1], when created using [3D3S] NI, create 9 Dice NI of a form that seems incompatible with this nomenclature, but we can nest our nomenclature just like we did the NIs to get a description: [3D3[3D3S]] NI.

\newpage

\section{Equivalency Between 5 Dice NI}

\begin{theorem}
The [5D:] and [5D:1] NI represent the same NI labelled differently.    

[5D:1] = [5D:] ADBEC

[5D:] = [5D:1] ACEBD
\end{theorem}

\noindent First, consider that all [5D:] and [5D:1] NI must be composed of overlapping [3D] NI. This can be drawn pictorially as follows: 

\begin{figure}[h!]
\includegraphics[width=\textwidth]{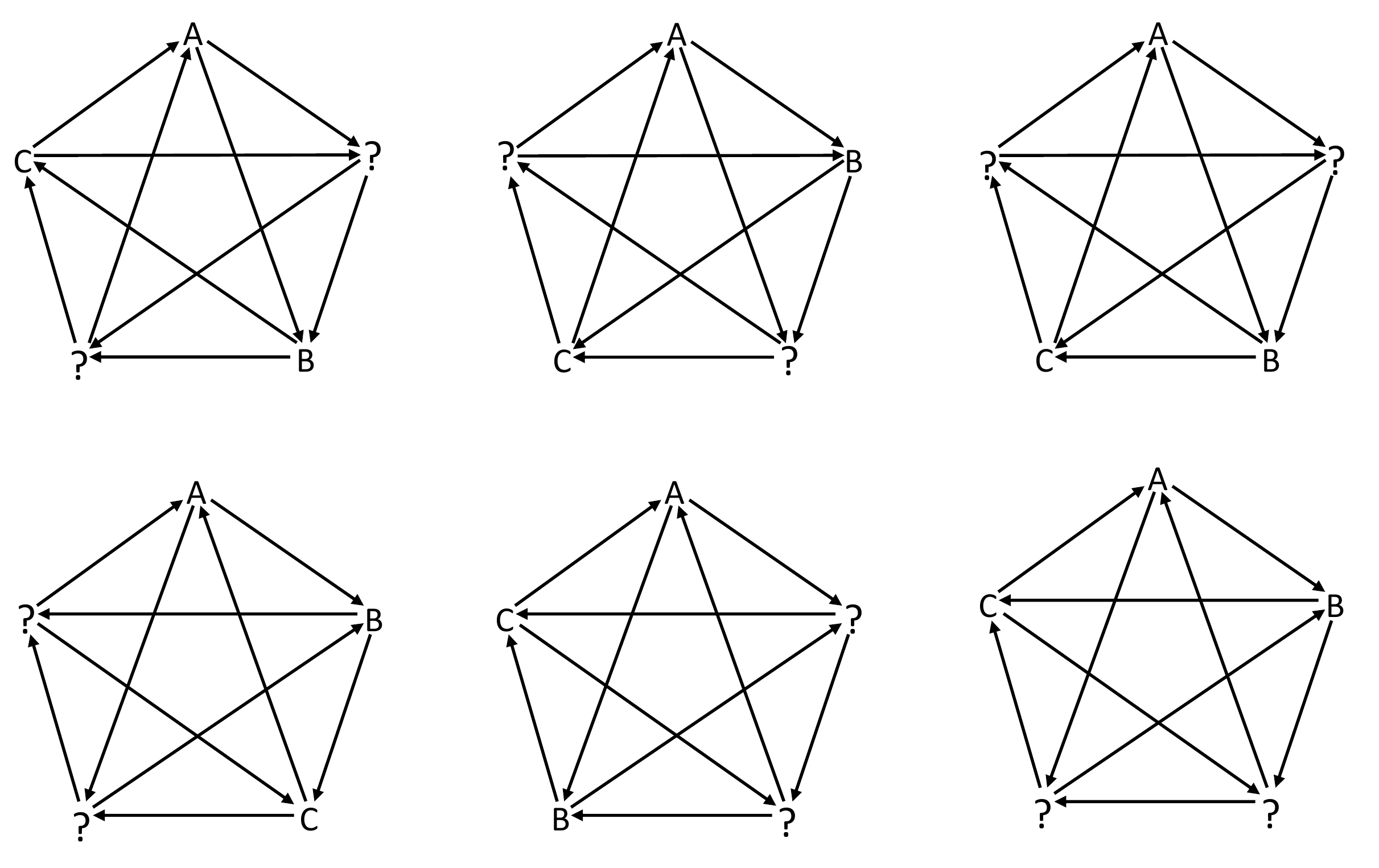}
\end{figure}

\noindent In this image, the 2 dice added to the identity are represented by the vertices marked with a \textquotesingle \textquotesingle ?\textquotesingle \textquotesingle. Now let us assume that our two added dice will follow the win chain imposed by this image:

\begin{figure}[h!]
\includegraphics[width=\textwidth]{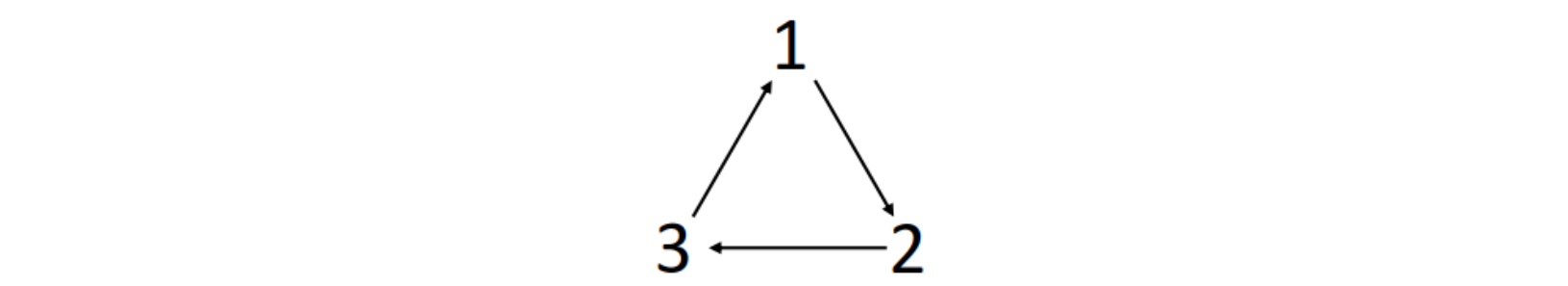}
\end{figure}

This assumption is useful because the pentagrams above predict who beats who, and if we know the win chain of our numbered dice, then, when we place one numbered die on the pentagram, we know what the other \textquotesingle \textquotesingle ?\textquotesingle \textquotesingle \  die must be. Now, knowing that Die A must beat one of the dies, we will arbitrarily define it to always beat die 1. With this definition, we can extrapolate as follows. The left most NI receive die 2 because die 2 is beat by die 1. The right most NI receive die 3 because die 3 beats die 1. The center column of NI receives a die 2 and a die 3 because that\textquotesingle s the only combination left.

\begin{figure}[h!]
\includegraphics[width=\textwidth]{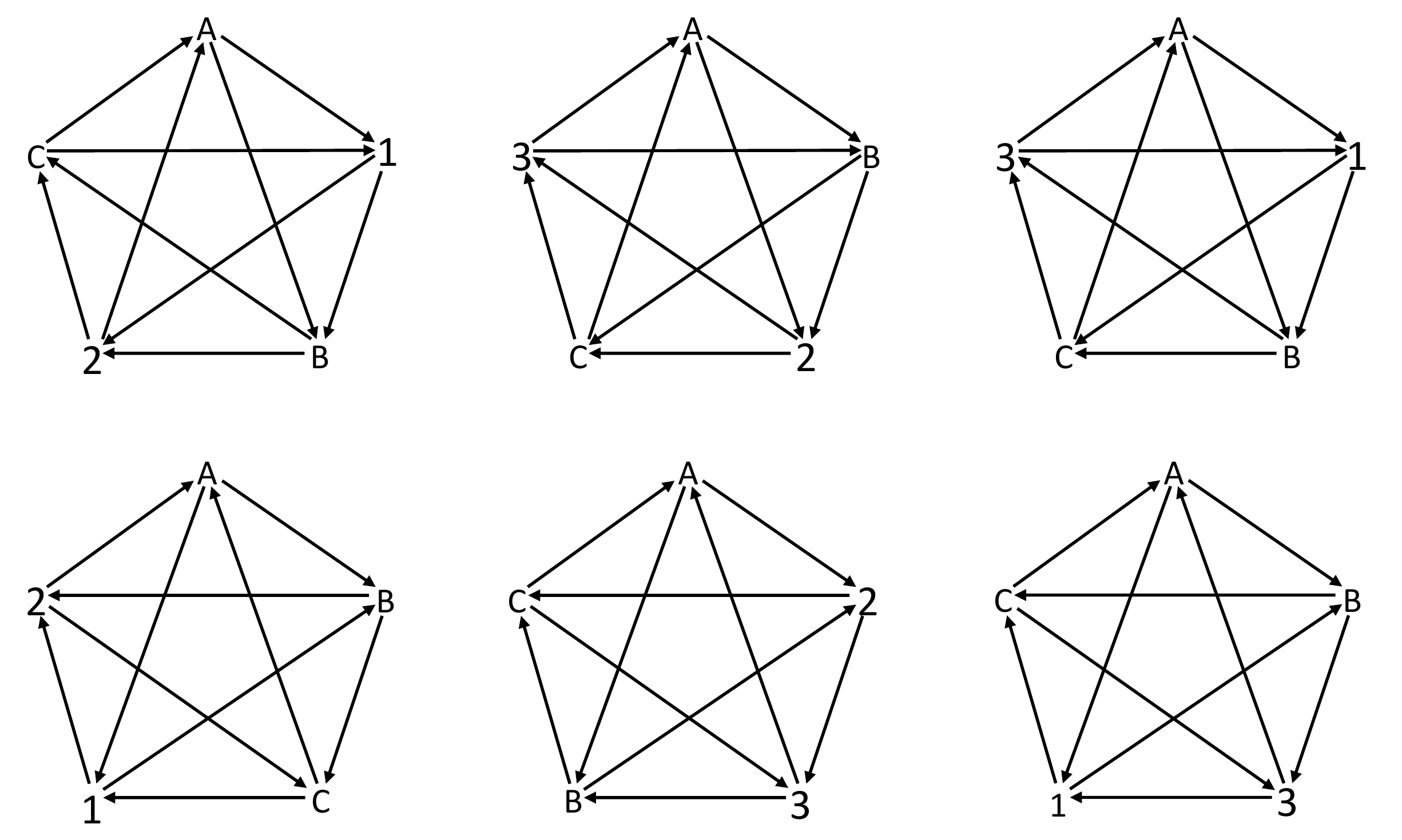}
\end{figure}

\newpage

\noindent When we look at the how the dice in these patterns perform against each other, we find that the [5D:] NI and [5D:1] NI produce the same patterns of wins, regardless of which win chain step pattern was used.

\begin{figure}[h!]
\includegraphics[width=\textwidth]{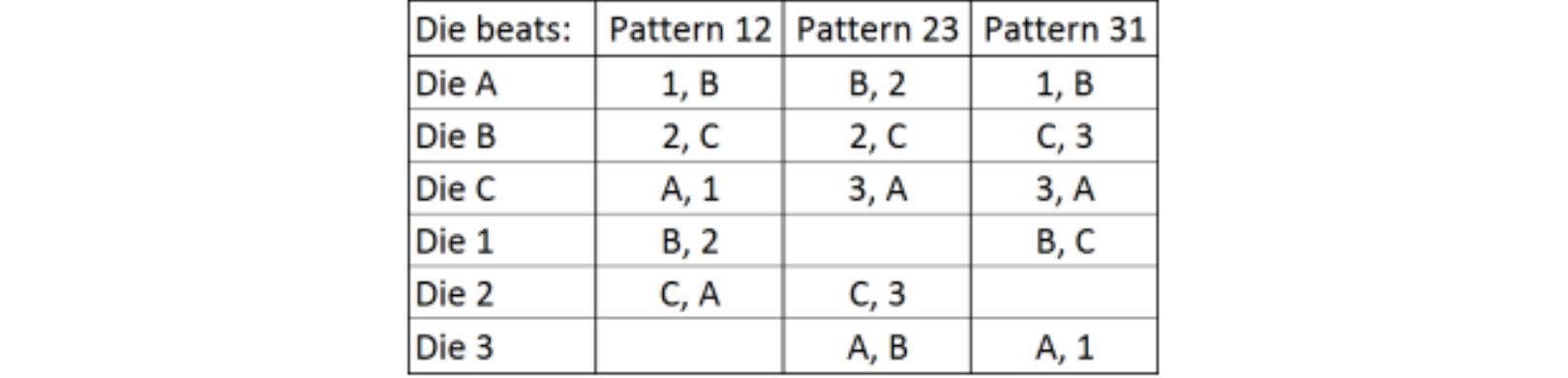}
\end{figure}

\noindent With this, we have proven that the [5D:1] NI are equivalent forms of the [5D:] NI and that the perceived difference can be removed by rearranging the primary win chain.

   [5D:] A1B2C = [5D:1] ABC12 

   [5D:] AB2C3 = [5D:1] A23BC

   [5D:] A1BC3 = [5D:1] AB31C

\

\noindent Stated generally:

   [5D:1] = [5D:] ADBEC

   [5D:] = [5D:1] ACEBD

\pagebreak

\noindent As such, a [5D:] NI will display a [5D:1] complete graph, if the [5D:] NI\textquotesingle s dies are reassigned as follows:

   Die B becomes Die D.

   Die C becomes Die B.

   Die D becomes Die E.

   Die E becomes Die C.

\ 

\noindent And a [5D:1] NI will display a [5D:] complete graph, if the [5D:1] NI\textquotesingle s dies are reassigned as follows:

  Die B becomes Die C.

  Die C becomes Die E.

  Die D becomes Die B.

  Die E becomes Die D.

\

\section{Adding Two Dice to a [3D] NI} 		
\
\noindent A [5D:] NI has five [3D] \textquotesingle \textquotesingle composing\textquotesingle \textquotesingle \  NI as follows:
\
\begin{figure}[h!]
\includegraphics[width=\textwidth]{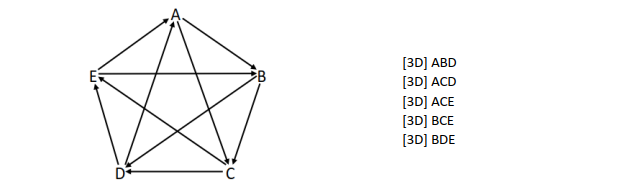}
\end{figure}

\noindent Definition: A NI is a \textquotesingle \textquotesingle composing\textquotesingle \textquotesingle \  NI when it is used in the composition of another NI.

\

The [3D] composing NI that involve the [5D:] die A, by default, have the die with the lowest face in the [5D:] NI, as such we can say that these [3D] NI will be of the following form:

[3D] ABD 

[3D] ACD 

[3D] ACE 

\

\noindent Definition: \textquotesingle \textquotesingle Anchored\textquotesingle \textquotesingle \ composing NI are composing NI that involve the Die A of the NI being composed.  

\

The remaining two [3D] NIs that compose the [5D:] NI lack a [5D:] NI die A to, in effect, \textquotesingle \textquotesingle anchor them\textquotesingle \textquotesingle \ to our NI list\textquotesingle s assumption of which die of a set has the face with the lowest value. As such, the remaining NIs can take multiple forms and can be referred to as \textquotesingle \textquotesingle unanchored\textquotesingle \textquotesingle .

[3D] BCE or [3D] CEB or [3D] EBC 

[3D] BDE or [3D] DEB or [3D] EBD 

\

Using this knowledge, we can create a [5D:] NI from [3D] NI, as seen in the following example. This example was originally solved in reverse using a [5D3S:] NI that was generated by modifying an almost perfectly nontransitive identity produced by the formula described in [2]. As such, this example proves that the theory described above works, but it does not provide much guidance in how to apply said theory with the highest chance of success.

\

Step 1: Assign a [3D] NI to all the \textquotesingle \textquotesingle anchored\textquotesingle \textquotesingle \ [3D] NIs of our theoretical [5D:] NI and reassign the letters as required. In this example, we will use the [3D] NI of A\textless C\textless B\textless C\textless B\textless A\textless B\textless A\textless C for all 3 anchored NI, this choice is not necessarily a requirement, but it will certainly make finding overlap between the [3D] NI much easier.

\begin{figure}[h!]
\includegraphics[width=\textwidth]{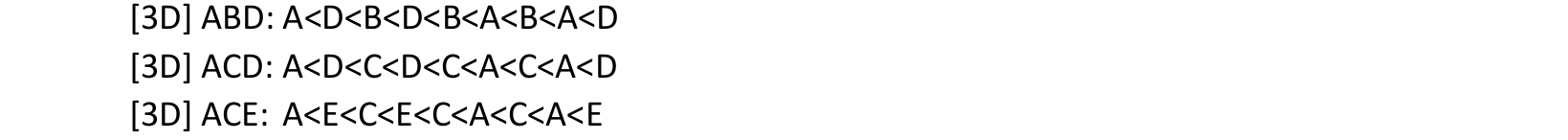}
\end{figure}

Step 2: Anchor our \textquotesingle \textquotesingle unanchored\textquotesingle \textquotesingle \ [3D] NIs and assign them all a [3D] NI. In this example, we will anchor [3D] BCE to [3D] BCE and [3D] BDE will be anchored to [3D] DEB. These [3D] NIs will then be assigned the following NIs, respectively, A\textless C\textless B\textless C\textless B\textless A\textless B\textless A\textless C and A\textless C\textless B\textless B\textless A\textless C\textless C\textless B\textless A.

\begin{figure}[h!]
\includegraphics[width=\textwidth]{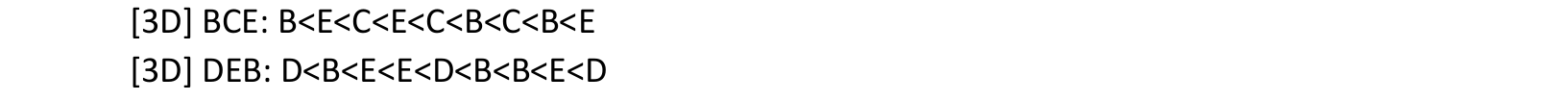}
\end{figure}

Step 3: Overlap the chosen NIs to generate the [5D:] NI with the appropriate number of sides. Note, each die face is overlapped an equivalent number of times.

\begin{figure}[h!]
\includegraphics[width=\textwidth]{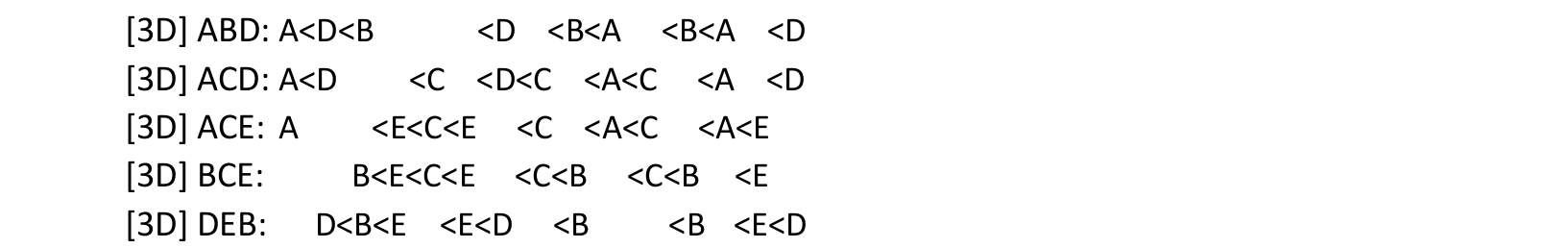}
\end{figure}

\noindent This produces the [5D3S:] NI: A\textless D\textless B\textless E\textless C\textless E\textless D\textless C\textless B\textless A\textless C\textless B\textless A\textless E\textless D.

\newpage

\section{Adding Two Dice to a [5D] NI} 
A [7D::] NI has seven [5D:] NI as follows:

\begin{figure}[h!]
\includegraphics[width=\textwidth]{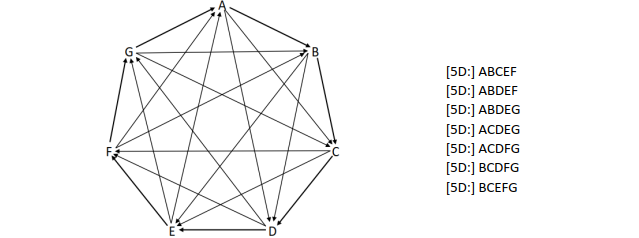}
\end{figure}

As before, we have a quantity of anchored composing NI equivalent to the number of dice involved in the composing NI, assuming the composing NI involve 2 less dice than our NI being generated. Also, as before, our new NI requires the use of two unanchored composing NI. 

Please note, these composing NI are correct for a [7D] NI of the form [7D::]. Given that the [5D] NI can appear different while being the same, we can assume that different win chain step patterns of high dice count NI can appear different while being the same.     

\begin{remark} 
The [9D] NI generated with the Identity Dice Multiplication method using [3D] NI, do not seem to be able to generate [5D] NI via any pattern of dice removal.
\end{remark}

\section{Patterns in the Frequency of Nontransitive Identities}

\noindent To iterate over identities, numbers were chosen to represent each operator and dice face combination. This was done as follows:

\begin{figure}[h!]
\includegraphics[width=\textwidth]{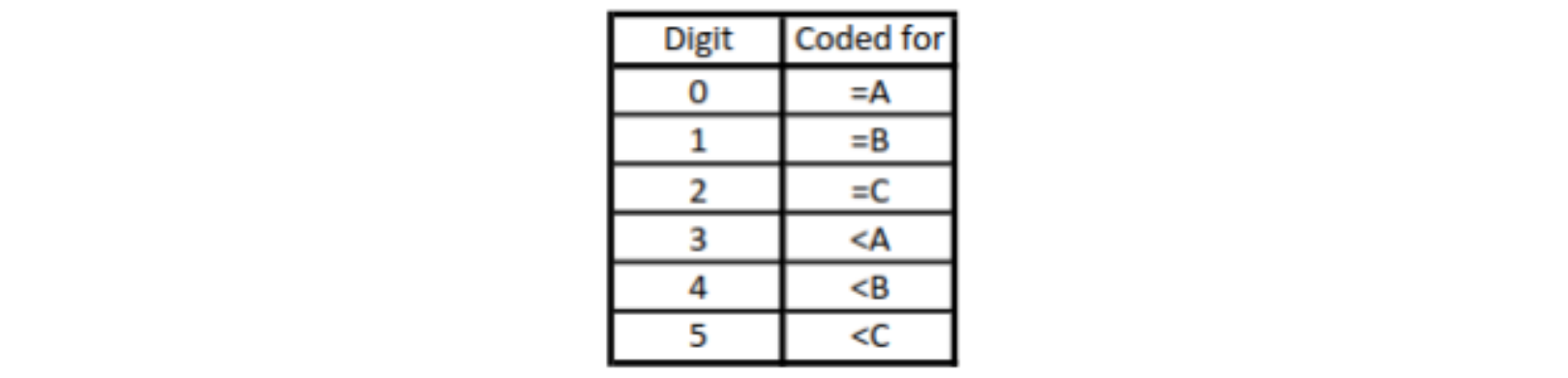}
\end{figure}

As such, for the 3 dice identities, we essentially iterating over a base-6 number. The script would start with the lowest possible \textquotesingle \textquotesingle viable\textquotesingle \textquotesingle \ identity and continue iterating that base-6 number until it was complete. For the [3D3S] NI, that lowest number was \textquotesingle \textquotesingle 00111222\textquotesingle \textquotesingle ,  of which represented the \textquotesingle \textquotesingle viable\textquotesingle \textquotesingle \  identity \textquotesingle \textquotesingle A=A=A=B=B=B=C=C=C\textquotesingle \textquotesingle . As die A is defined to have the lowest number die face, it\textquotesingle s presence at the beginning of the identity is assumed and therefore there is no reason to represent it in our iterated number. 

\begin{remark} 
A viable identity may or may not be nontransitive. A viable identity is only defined as having met the requirement of having the appropriate number of die faces for each die.
\end{remark}

Having started with \textquotesingle \textquotesingle 00111222\textquotesingle \textquotesingle ,  the next number to be tested would be \textquotesingle \textquotesingle 00111223\textquotesingle \textquotesingle ,  this number would fail the viability check because it has too many faces for die A (“A=A=A=B=B=B=C=C\textless A”), as such it would be skipped. The next number, \textquotesingle \textquotesingle 00111224\textquotesingle \textquotesingle ,  would also be skipped for having too many sides on dice B. The process of adding 1 to our base-6 iterated number, checking whether it was viable identity and if it was, whether it was also a nontransitive identity, would continue until our iterated number reached \textquotesingle \textquotesingle 55500000\textquotesingle \textquotesingle \  because, at this point, all further viable identities produced would feature a Die C that lost to Die A, and this result would violate the defined nontransitive relationship we were seeking. 

With this, I gained two lists of numbers. The first list coded for the viable identities and the second list coded for the nontransitive identities. Starting with the unproven assumption that the \textquotesingle \textquotesingle nonviable\textquotesingle \textquotesingle \ identities were essentially garbage numbers, I chose to set an index to the viable identities. As such, the first viable identity \textquotesingle \textquotesingle 00111222\textquotesingle \textquotesingle \ received an index of \textquotesingle \textquotesingle 0\textquotesingle \textquotesingle \ and the next viable identity received an index of 1 and so on. These \textquotesingle \textquotesingle viable identity indexes\textquotesingle \textquotesingle \ allowed me to assign the nontransitive identities a spot along a spectrum that was defined as being viable. I then chose to look for a pattern in the difference between the successive viable identity indexes of each NI.

\begin{remark} 
As far as NI list stringency, there are two possible goals. Any published list of NI should be listed such that the list is irreducible, i.e. no NI in the list is essentially a duplicate of another NI in the list. For example, any NI involving a sequence of the form \textquotesingle \textquotesingle \textless N\textless N\textquotesingle \textquotesingle  where N refers to a face of any one die, is essentially the same as one with the sequence \textquotesingle \textquotesingle \textless N=N\textquotesingle \textquotesingle ,  as such, all viable identities with a sequence of the form \textquotesingle \textquotesingle \textless N\textless N\textquotesingle \textquotesingle \ can be skipped for the purposes of a published list of NI because the same NI will be found and tested in the viable identity that features the \textquotesingle \textquotesingle \textless N=N\textquotesingle \textquotesingle \  sequence instead. The same situation exists for sequences that involve equivalent die faces. These faces can exist in a viable nontransitive identity in any order, but allowing any order produces many duplicative NI. Setting the requirement that NI feature equivalent die faces in alphabetical order restricts us to one allowed pattern for describing equivalent die faces and thus removes duplicative NI for this situation. This requirement is great for producing an irreducible list of NI, but duplicative NI are no less nontransitive than any other NI. For this reason, if our goal is to find a pattern in the frequency of NI, then we must consider the duplicative NI in addition to those of the irreducible list. This article will consider the patterns of both the irreducible and duplicative lists of [3D3S] NI and [3D4S] NI, the tables will be labelled to make sure this distinction is clear. Please note, the \textquotesingle \textquotesingle viable\textquotesingle \textquotesingle \  identities for the irreducible list of [3D3S] NI were held to the same requirements as the NI, as such their list of viable identities is a \textquotesingle \textquotesingle irreducible\textquotesingle \textquotesingle \ list of viable identities.
\end{remark}

\begin{figure}[h!]
\includegraphics[width=\textwidth]{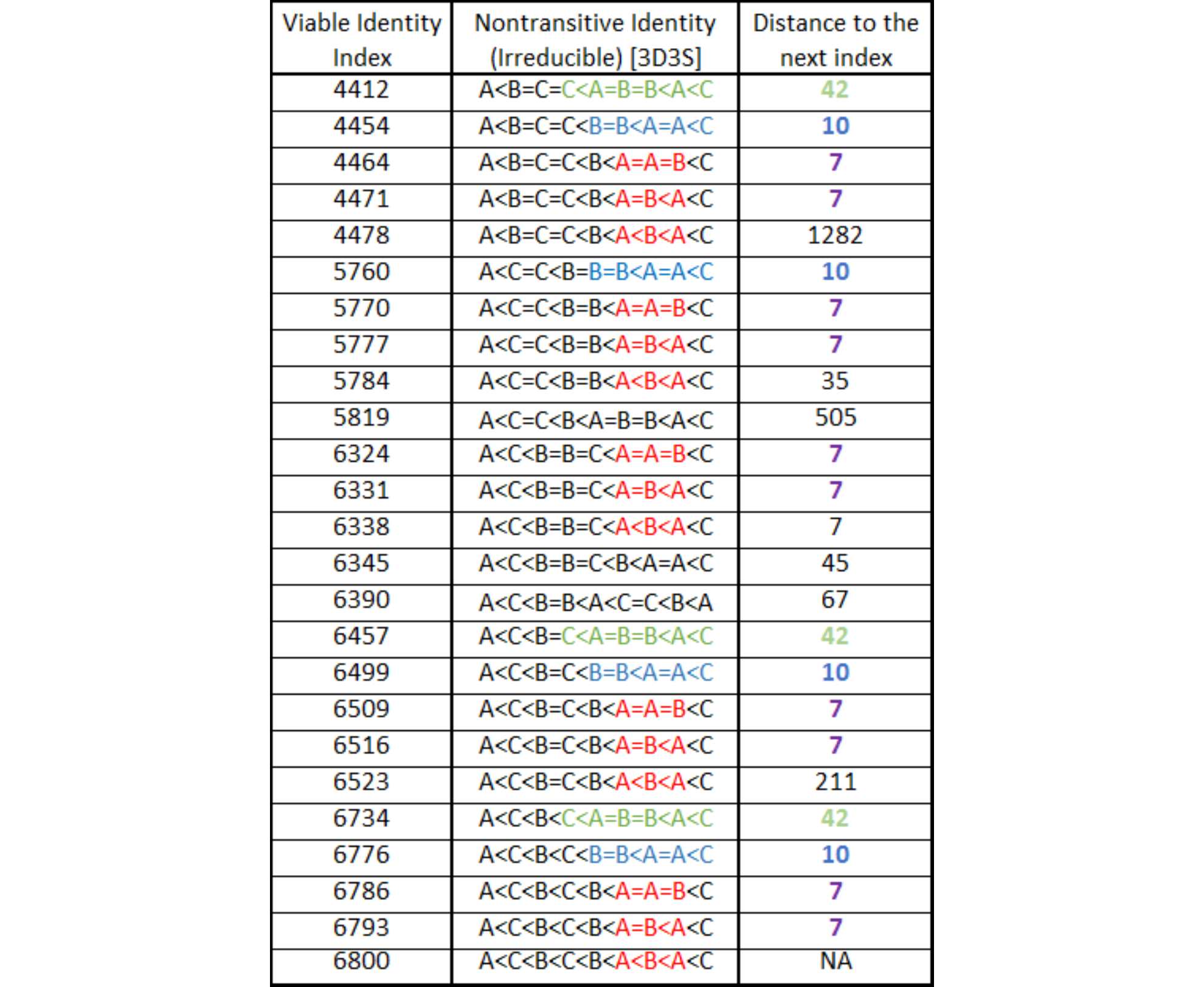}
\end{figure}

\newpage

The patterns in the distance between successive NI viable identity indexes are rather easy to see in this chart and they are especially clear in the color copy version, as the relevant distances and portions of the NI for each pattern have been color-coded. The pattern \textquotesingle \textquotesingle 42, 10, 7, 7\textquotesingle \textquotesingle \ repeats 5 times, 2 of which are partial repeats. Closer inspection notes that these partial repeats are due a \textquotesingle \textquotesingle full\textquotesingle \textquotesingle \  repetition requiring the production of a duplicitous NI, of which we have banned in this set of NI. We can also note that these patterns correlate to specific changes in portions of the NI from another NI. As such, we know that these distance patterns correlate to our NI\textquotesingle s morphisms. 

\begin{remark} 
The remaining charts are a bit large for direct publication, they have been included in this article\textquotesingle s ancillary files in the \textquotesingle \textquotesingle Frequency of NI tables.xlsx \textquotesingle \textquotesingle \ document.
\end{remark}

The pattern in the duplicitous [3D3S] NI chart is also rather immediately apparent. There are 13 repetitions of the pattern \textquotesingle \textquotesingle 8, 4, 24, 12, 4, 8, 4\textquotesingle \textquotesingle ,  of which, 10 are preceded by a \textquotesingle \textquotesingle 32\textquotesingle \textquotesingle . Additionally, there are 8 repetitions of the pattern \textquotesingle \textquotesingle 168, 40, 8, 4\textquotesingle \textquotesingle ,  of which, 4 are preceded by a \textquotesingle \textquotesingle 164\textquotesingle \textquotesingle . 

\begin{remark} 
At this point, it feels as if I\textquotesingle m just throwing out numbers like they were going out of style. As such, I want to take a moment to note that there is a value behind understanding these numbers, that value is that if we can identify the pattern behind this set, then we can generalize that pattern to identify the NI of sets that would otherwise take an eternity to compute.
\end{remark}

The pattern for the irreducible list of [3D4S] NIs are less friendly. Using this method, the patterns seen only repeats once and part of the pattern can be interrupted with NI not predicted by the first or second repetition of the pattern. For example, the first pattern in the set is \textquotesingle \textquotesingle 2, 4, 1, 41, 2, 1, 2\textquotesingle \textquotesingle ,  when this pattern occurs again, it occurs as \textquotesingle \textquotesingle 2, 4, 1, 18, 2, 21, 2, 1, 2\textquotesingle \textquotesingle . Note, 18+2+21 = 41. As such, the pattern did hold, it just didn\textquotesingle t predict the NI that existed at +18 and +20 in the irreducible list of [3D4S] NIs.

\

The first pattern in the duplicative [3D4S] NI list is \textquotesingle \textquotesingle 3, 1, 1, 42, 1, 1, 38, 1, 1, 7, 3, 1, 1\textquotesingle \textquotesingle ,  it repeats 13 times in the first 170 NI, for 3 of these repeats, the \textquotesingle \textquotesingle 3, 1, 1, 42\textquotesingle \textquotesingle  are missing. The second series of patterns involves \textquotesingle \textquotesingle 1, 1, 1, 1, 1, 1, 1, 1, 1, 38, 1, 1, 1, 1, 4, 32, 1, 1, 1, 1, 2, 3, 1, 1, 1, 1, 1, 1, 1, 1, 1, 6, 6\textquotesingle \textquotesingle \  and \textquotesingle \textquotesingle 8, 32, 6, 18, 6\textquotesingle \textquotesingle ,  of which repeat 12 and 11 times respectively. Between each repetition of the pattern is either nothing or a large number. I checked the first 3100 distances by hand and these types of patterns keep occurring. I wrote a script to examine the patterns in a slightly different manner, a summary of these patterns is available in this article\textquotesingle s ancillary files. 

\

The existence of these patterns implies that the patterns in nontransitive identities exist on a smaller level than that of the NI themselves. Understanding the nature of these morphisms will help immensely in our search for more complicated NI. 

\section{Future Work}

Future work could be directed in many directions. Identifying the [3D6S] NI would be especially interesting as it could show which identities become possible with 6 sides that are not possible based on operations that involve [3D3S] NI. The list could also be used to identify other operations that can be performed with NI.

\

The identities for even numbers of dice should, in my opinion, be skipped as they will all simply be identities for odd numbers of dice with one die removed. Writing a program that identifies [5D] NI based on the [3D] NI would likely be significantly faster than the script that I used to generate my [5D3S] NI list. Identification of identical forms of [7D] NI and [9D] NI would also be useful. It\textquotesingle s possible that the second is already known, as it relates to complete graphs but I wasn\textquotesingle t able to find anything in my search of the literature. If it does exist, then the translation of this research into the nomenclature of Nontransitive Identities would be useful.

\

A method for expanding NI that I did not mention in this paper involves solving multiple NI individually and then mapping those individual solutions to one set of dice. The idea behind this method being that the NI where being solved \textquotesingle \textquotesingle in parallel\textquotesingle \textquotesingle ,  instead of \textquotesingle \textquotesingle in series\textquotesingle \textquotesingle  as was the case with the Identity Addition method. I found that this method worked for a subset of NIs, but I was unable to characterize what factors led to it working or not working. I also tested a modification of the \textquotesingle \textquotesingle Parallel Identity Addition\textquotesingle \textquotesingle \  method, where a greater increase between subsequent die faces was used for subsequent \textquotesingle \textquotesingle parallel\textquotesingle \textquotesingle \ NI, what was of interest here was that the subset of NI that this \textquotesingle \textquotesingle Staggered Parallel Identity Addition\textquotesingle \textquotesingle \ method worked for seemed to be a subset of the NI that the \textquotesingle \textquotesingle Parallel Identity Addition\textquotesingle \textquotesingle \  method worked for.

\

The morphisms of nontransitive identities are another area for future work. These morphisms could be characterized and further developed in terms of higher-dimensional categories.

\section{Concluding Remarks}

I would like to thank B. Haran of the youtube channel Numberphile as well as many of the mathematicians that he\textquotesingle s interviewed (M. Parker, J. Grime, and T. Tokieda) for inspiring this research. I would also like to thank the users of Stackoverflow\textquotesingle s python section. I did not know how to program in python prior to this project and I would not have been able to write my scripts without the time they spent answering other people\textquotesingle s questions. Also, I must apologize in advance if I take a while to respond to emails, I start my third year of medical school on May 8th.

\section{Correspondence}

\noindent logandanielson@gmail.com

\section{References}

\noindent [1] Martin Gardner, \textquotesingle \textquotesingle The Paradox of the Nontransitive Dice and the Elusive Principle of Indiﬀerence\textquotesingle \textquotesingle ,  Scientiﬁc American (1970), 223(6) pp. 110–114.

\noindent [2] J.W. Moon, L. Moser, Generating oriented graphs by means of team comparisons, Pacific Journal of Mathematics (1967), 21(3) pp. 531–535

\end{document}